\documentclass{alf-j-l}
\usepackage{verbatim} 
\usepackage{url} 
\usepackage{AlfPreambleL}  

\theoremstyle{plain}  
\newtheorem{theorem}{Theorem}

\newtheorem{corollary}[theorem]{Corollary}

\newtheorem*{theorem*}{Theorem}

\theoremstyle{definition}

\newtheorem{example*}{Example}

\theoremstyle{remark}

\newtheorem*{Remarks*}{Remarks} 
\newtheorem*{Remark*}{Remark}

\PII{Centre for Number Theory Research, 1 Bimbil Place, Killara,
Australia 2071}%

\copyrightinfo{\number\year}{Francesco Pappalardi and Alfred J van
der Poorten}

\begin{document}

\def\currentvolume{173}
\def\currentissue{Draft of\ }
\def\paperdate{\today}
\def\ISSN{}
\pagespan{69}{\pageref{page:lastpage}}   

\def\ceNTRelogo{\relax}

\title {Pseudo-Elliptic
Integrals,\rlap{\smash{\quad\qquad\ceNTRelogo}}\\ Units, and 
Torsion\\\phantom{nothing}} 

\author{Francesco Pappalardi}
\address{Dipartimento di Matematica, 
Universit\`a degli Studi Roma Tre\newline 
L.go San Leonardo Murialdo, 1 
c.a.p. I-00146 Roma, Italy}
\email{pappa@mat.uniroma3.it  (Francesco Pappalardi)}

\author{Alfred J. van der Poorten}
\address{Centre for Number Theory Research\newline 
1 Bimbil Place, Killara,
Sydney NSW 2071, Australia}
\email{alf@math.mq.edu.au (Alf van der Poorten)}

\thanks{This paper was constructed during a visit by the second
author to Italy supported in part by a GNSAGA INDAM grant; his work
is also supported by a grant from the Australian Research Council.}



\date{\Now.}


\begin{abstract} We remark on pseudo-elliptic integrals and on
exceptional function fields, namely function fields defined over an
infinite base field but nonetheless containing non-trivial units.
Our emphasis is on some elementary criteria that must be satisfied
by a squarefree \poly\ whose square root generates a
quadratic function field with non-trivial unit. We detail the genus~$1$ case.
 \end{abstract}

\maketitle
\pagestyle{myheadings}\markboth{{\headlinefont Francesco
Pappalardi and Alf van der Poorten}}{{\headlinefont Pseudo-elliptic
integrals}}

\section{Pseudo-Elliptic Integrals}
\noindent  The surprising integral
\begin{multline*}
\int\frac{6x\,dx}{\sqrt{x^4+4x^3-6x^2+4x+1}}=
\log\Bigl(x^6 + 12x^5 + {45}x^4+ 44x^3 - {33}x^2 + 43\\+
(x^4 + 10x^3 + 30x^2 + 22x - 11)\sqrt{x^4+4x^3-6x^2+4x+1}\,\Bigr)
\end{multline*}
is a nice example of a class of pseudo-elliptic integrals
\begin{equation}\label{eq:pseudo}
\int\frac{f(x)dx}{\sqrt{D(x)}}=\log\(a(x)+b(x)\sqrt{D(x)}\).
\end{equation}
Here we take $D$ to be a monic \poly\ defined over
$\Q$, of even degree $2g+2$, and  not the square of a \poly;  $f$,
$a$, and $b$ denote appropriate
\poly s. We suppose $a$ to be nonzero, say of degree $m$ at least
$g+1$. We will see that necessarily $\deg b=m-g-1$, that $\deg f=g$,
and that
$f$ has leading coefficient $m$. In our example, $m=6$ and $g=1$.

Plainly, if \eqref{eq:pseudo} holds then it remains true with
$\sqrt D$ replaced by its conjugate $-\sqrt D$. Adding the two
conjugate identities we see that
\begin{equation}\label{eq:unit1}
\int 0\,dx=\log\(a^{2}-Db^{2}\).
\end{equation}Thus $a^{2}-Db^{2}$ is some constant $k$, and must be
nonzero because $D$ is not a square. In other words, $u=a+b\sqrt D$
is a nontrivial unit in the function field $\Q(x,\sqrt {D})$;
 and $\deg a=m$ implies $\deg b=m-g-1$ is immediate.
 
Differentiating \eqref{eq:unit1} yields $2aa'-2bb'D-b^{2}D'=0$.
Hence $b\Div aa'$, and since $a$ and $b$ must be relatively prime 
because $u$ is a unit, it follows that $b\Div a'$. Set $f=a'/b$,
noting that indeed $\deg f=g$ and that $f$ has leading coefficient
$m$ because $a$ and $b$ must have the same leading
coefficient.\footnote{That common coefficient is $1$ without loss of
generality since we may freely choose the constant produced by the
indefinite integration.}
 
Moreover,
  \begin{equation*}
u'=a'+b'\sqrt D+bD'/2\sqrt D
=a'+(2bb'D+b^{2}D')/2b\sqrt D=a'+aa'/b\sqrt D\,.
\end{equation*}
So, remarkably, $u'=f(b\sqrt D +a)/\sqrt D=fu/\sqrt D$.  
  
Thus, to verify \eqref{eq:pseudo} it suffices to make the not
altogether obvious substitution
$u(x)=a+b\sqrt D$, of course given that $u$ is a unit of the order
$\Q[x,\sqrt {D}\,]$.

\begin{Remark*} The case $g=0$, say $D(x)=x^2+2vx+w$, is
useful for orienting oneself. Here $(x+v)+\sqrt D$ is a
unit, of norm $v^2-w$, and indeed
$$\int\frac{dx}{\sqrt{x^2+2vx+w}}=\arsinh\frac{x+v}{\sqrt{w-v^2}}
=\log\(x+v+\sqrt{x^2+2vx+w}\,\)\,.
$$
Notice that $\deg f=0$ and has leading coefficient $1$, as predicted.\end{Remark*}

\section{Units in Quadratic Extension Fields, and Torsion}

\subsection{Number fields} Let $N$ be a positive integer, not a
square, and set
$\omega=\sqrt N$. It is easy to apply the Dirichlet box principle to
prove that an order
$\Z[\omega]$ of a quadratic number field $\Q(\omega)$ contains
nontrivial units. 
Indeed, by that principle there are infinitely
many pairs of integers $(p,q)$ so that $|q\omega-p|\lt1/q$, whence
$|p^2-Nq^2|\lt 2\sqrt N+1$. It follows, again by the box principle, 
that there is an integer
$l$ with $0\lt |l|\lt 2\sqrt N+1$ so that the equation $p^2-Nq^2=l$
has infinitely many pairs $(p,q)$ and $(p',q')$ of solutions with
$p\=p'$ and $q\=q'\pmod l$. For each such distinct pair,
$al=pp'-Nqq'\!$,
$bl=pq'-p'q$, yields $a^2-Nb^2=1$.

\subsection{Function fields} Just so, in the function field case
already introduced, there are infinitely many pairs of \poly s
$p(x)$ and $q(x)$ so that
$\deg(q\sqrt D-p)\lt-\deg q$, whence $\deg(p^2-Dq^2)\le g$. But a
second application of the box principle fails when the base field,
$\Q$ in our introductory discussion, is infinite; because there are
then infinitely many distinct \poly s of bounded degree. In that
case, the existence of a nontrivial unit (thus, one not an element
of the base field) is unusual happenstance. Accordingly,
we say that a function field $\Q(x, \sqrt D)$  with a nontrivial
unit $a+b\sqrt D$ is an
\emph{exceptional} function field and we call $D$ an
\emph{exceptional} \poly.

\subsection{Torsion on the Jacobian of a hyperelliptic curve} A
slight change of viewpoint, emphasising the hyperelliptic curve
$\mathcal C: y^2=D(x)$, may clarify matters. A function $u=a+by$ is
a unit precisely if its divisor is supported only at infinity.
However, $\mathcal C$ has two points at infinity, say $O$ and $S$
(or $\infty_-$ and $\infty_+$ if one prefers) and so the divisor of
$u$ is some multiple, say $m(S-O)$, of the divisor $S-O$ at
infinity. Because $u$ is a function, this is to say that the class
of $S-O$ on the Jacobian of $\mathcal C$ is torsion of order $m$.
In the case
$\deg D=4$, so genus $g=1$ if $D$ is squarefree, we may take $O$ as
the zero of the elliptic curve $\mathcal C$ and report that the
point $S$ on $\mathcal C$ is torsion of order $m=\deg a$.

\section{Exceptional Quadratic Fields.}

\noindent It is appropriate to identify straightforward
properties of the squarefree \poly\ $D(x)=y^2$ sufficient or just
necessary  that the field $\Q(x,y)$ be exceptional.

Suppose therefore that $\Q(x,y)$ is exceptional, so that we have a
unit $u=a+by$ or, more helpfully, an identity $b^2D=a^2-k$ with
$a,b\in \Q[x]$ and $k\in\Q\setminus\{0\}$. It will be helpful to
set $k=c^2$. We note immediately that the two \poly s $a-c$ and
$a+c$, which are conjugate over $\Q$ if $k$ is not a square,  are
relatively prime.

We have $b^2D=(a-c)(a+c)$. Hence if $k$ is not a square, 
$b$ must factor in $\Q(c)[x]$ as a norm $d\conj d$, where the
overline
$\conj{\phantom{d}}$ denotes conjugation in the quadratic extension
$\Q(c)$, and $D$ factorises over $\Q(c)$ as the product of the
\poly\ $(a-c)/d^2$, and of its conjugate. In particular, $\deg
b=m-g-1$ must be even.

If, however, $k$ is a square in $\Q$ (thus, in particular,
always if
$\deg b=m-g-1$ is odd) then we seem to see only that $b$ must have a
factor
$d$ defined over $\Q$ so that both $2\deg d$ and $2m -(2g+2)-2\deg
d$ do not exceed $m=\deg (a-c)=\deg a$. That is, we have
$m-(2g+2)\le 2\deg d
\le m$.

\begin{theorem} \label{th:1} Set $y^2=D(x)$, with $D$
monic, squarefree, and of degree $2g+2$. Suppose the domain
$\Q[x,y]$ 
contains a unit of degree
$m\gt g$ and norm $k$.
\begin{itemize}
\item[(a)] If $m$ and $g$ have the same parity then $k=c^2$ is a
square.
\item[(b)] If $k=c^2$ is a square, there is a positive integer $s$
so that 
$D$ is a product of \poly s over $\Q$ of degrees $m-2s$ and
$2g+2+2s-m$. Thus $D$ is reducible over $\Q$ if 
$m$ is odd.
\item[(c)] If $k=c^2$ is not a square in $\Q$ then $D$ factorises
over $\Q(c)$ as a product of two polynomials conjugate over
$\Q(c)$, so each of degree $g+1$. 
\end{itemize}
\end{theorem}

\noindent Note that the compactly written assertion (b) includes the
possibility that $D$ is irreducible if $m$ is even, and (since both
the stated degrees must be nonnegative) that it implicitly entails
upper and lower bounds on the integer $s$. Assertion (c) implies
that the Galois group of
$D$ is restricted by $\#\Gal(D)\Div 2\((g+1)!\)^2$. Thus, if $g=1$
it is the dihedral group on four elements or one of its subgroups.
 
We observe also that the
statements of the theorem, which refer only to the
\poly\ $D$ and the torsion order~$m$, do not include all the
information that may be extracted from the remarks preceding the
proclamation of the theorem.

\begin{Remarks*} It should be no surprise that none of the
criteria of the theorem suffice to guarantee obtaining an
exceptional quadratic function field. We detail the case $g=1$ in
\S\ref{g1} at page~\pageref{g1} below.
 \end{Remarks*}

\section{Continued Fractions}
\subsection{Number fields} There is a well known algorithm in the
number field case yielding the fundamental unit of the order
$\Z[\sqrt N]$. As before set $\omega=\sqrt N$ and suppose
$A$ is the integer part of $\omega$. The zero-th step in the \cfe\
of
$\omega+A$  is
\begin{equation}\label{eq:linezero}
\omega+A =2A -(\conj\omega+A)\end{equation}
and a typical consequent step is of the shape
$$(\omega+P_h)/Q_h=a_h-(\conj\omega+P_{h+1})/Q_h\,; \quad\text{in
brief}\quad \omega_h=a_h-\conj\rho_h\,. 
$$
Thus $P_h+P_{h+1}+(\omega+\conj\omega)=a_hQ_h$, and because the next
complete quotient $\omega_{h+1}$ is the
reciprocal of the remainder $-\conj\rho_h$ we must
also have
$(\omega+P_{h+1})(\conj\omega+P_{h+1})=-Q_hQ_{h+1}$. In particular,
certainly $Q_{h+1}$ divides the norm
$(\omega+P_{h+1})(\conj\omega+P_{h+1})$.

Here the $P_h$ and $Q_h$ are integers, and it is readily shown they
all satisfy
\begin{equation}\label{eq:bounds}
0\lt 2P_h+(\omega+\conj\omega)\lt \omega-\conj\omega\,,\qquad
0\lt Q_h \lt \omega-\conj\omega
\end{equation}
proving, by the box principle, that the \cfe\ of $\omega$ is
periodic. Moreover, one notices that always both
\begin{equation}\label{eq:reduced}
\text{$\omega_h\gt1$  while
$-1\lt \conj \omega_h\lt 0$, and $\rho_h\gt1$  while
$-1\lt \conj \rho_h\lt 0$}.
\end{equation} 
It follows that
conjugation of the \cf\ tableau, replacing  
$$\omega_h=a_h-\conj\rho_h \quad\text{by}\quad
\rho_h=a_h-\conj\omega_h\,,$$
again gives a \cfe\ --- in particular, $a_h$ which began life as the
integer part of $\omega_h$, also is the integer part of $\rho_h$ ---
reversing the order of the lines of the original \ex. Because
line zero \eqref{eq:linezero} is symmetric it occurs in the \ex\ of
$\rho_h$, and because the \ex\ of $\omega+A$ is
periodic it follows that it is in fact pure periodic, moreover with
a symmetry: if the period length is $r$ then the word
$\cCF{a_1\\a_2\\\ldots\\a_{r-1}}$ must be a palindrome. 

One obtains the fundamental unit $a+b\omega$ by computing the
convergent
\begin{equation}\label{eq:unit}
\CF{A\\a_1\\a_2\\\ldots\\a_{r-1}}=a/b\,.
\end{equation}

\subsection{Function fields} Mutatis mutandis, the function field
argument is identical. We set
$y^2=D(x)$ as before. Plainly we may write $D$ as $D=A^2+R$, where
$\deg A=g+1$ and $\deg R\lt g$; then $A$ is the \poly\ part of the
Laurent series $y\in\Q((x^{-1}))$. We expand $y+A$ in complete
analogy with the numerical case, but now selecting the \pq s $a_h$
as the \poly\ part of the respective complete quotients
$y_h:=(y+P_h)/Q_h$. 
The bounds \eqref{eq:bounds} become
\begin{equation}\label{eq:boundsf}\deg P_h=g+1 \quad\text{and}\quad
\deg Q_h\le g
\tag{$\ref{eq:bounds}^\prime$}\end{equation}
and of course do not guarantee periodicity, because the base field
$\Q$ is infinite. The conditions \eqref{eq:reduced} for reduction
turn into
\begin{multline}\label{eq:reducedf}
\deg (y+P_h)\gt \deg Q_h \quad\text{but}\quad \deg (\conj y+P_h)\lt
\deg Q_h\\\quad\text{and}\quad  \deg (y+P_{h+1})\gt \deg Q_h
\quad\text{but}\quad \deg (\conj y+P_{h+1})\lt \deg Q_h\,.
\tag{$\ref{eq:reduced}^\prime$}\end{multline}
As in the number field case, conjugation reverses the \cf\ tableaux.
Thus, \emph{if} the \ex\ of $y+A$ happens to be periodic then it has
the symmetries of the number field case and the \cfe\ yields a unit
of norm~$1$, given by the convergent \eqref{eq:unit}.

\subsection{Quasi-periodicity}\label{ss:quasi} Suppose now that $D$
is exceptional in that the function field $\Q(x, y)$ contains a
unit
$u$, of norm
$-\kappa$. By general principles that entails that some $Q_i$ is $\pm
\kappa$, say
$Q_r=\kappa$ with $r$ odd. That is, line~$r$ of the \cfe\ of $y+A$
is
\begin{equation}\label{eq:liner}
y_r:=(y+A)/\kappa=2A/\kappa-(\conj
y+A)/\kappa\,;\tag*{line~$r$:}\end{equation}  here we have used
\eqref{eq:reducedf} to deduce that necessarily
$P_r=P_{r+1}=A$. We recall that
\begin{equation}\label{eq:line0}
 y+A=2A-(\conj y+A)\,.\tag*{line~$0$:}\end{equation}
By conjugation of the $(r+1)$-line tableau showing that $y+A$ is
quasi-periodic we see immediately also that  
\begin{equation}\label{eq:line2r}
y_{2r}:=y+A=2A-(\conj y+A)\,,\tag*{line~$2r$:}\end{equation}
so that in any case if $y+A$ has a quasi-periodic \cfe\ then it
is periodic of period twice the quasi-period. This is
a result of Tom~Berry \cite{Be}; it applies to arbitrary quadratic
irrational functions whose trace%
\footnote{If $y$ has trace $t$, rather than zero trace, replace line
zero of the \ex\ by $y+A-t=2A-t-(\ov y+A-t)$ and so on in the story
just told. To be able to do that $t$ should of course be `integral',
that is, a \poly.}
is a \poly. Other
elements
$(y+P)/Q$ of
$\Q(x,y)$, with $Q$ dividing the norm   $(y+P)(\conj y+P)$, may
be honest-to-goodness quasi-periodic, that is, not also periodic. 

\label{page:quasiperiodic}

Further, if $\kappa\ne-1$ then $r$ \emph{must} be odd. To see that,
notice the identity
$$B\CF{Ca_0\\Ba_1\\Ca_2\\Ba_3\\\ldots}=C\CF{Ba_0\\Ca_1\\Ba_2\\Ca_3\\\ldots},
$$
reminding one how one multiplies a \cfe\ by some
quantity; this cute formulation of the multiplication rule is due
to Wolfgang Schmidt \cite{Schm}. The `twisted symmetry' occasioned
by division by $\kappa$, equivalent to the existence of a
non-trivial quasi-period, is noted by Christian Friesen \cite{Fr}. 

In summary, if quasi-periodic it is periodic, and the \cfe\ of
$y=\sqrt{D}$ has the symmetries of  the more familiar number field
case, as well as  twisted symmetries occasioned by a nontrivial
$\kappa$.

\begin{Remarks*} The conclusion just stated is surely
well known. Certainly it is asserted by Adams and
Razar \cite{AR}, but without the couple of lines of argument we add
here. The second of us is indebted to notes of Ethan Street
\cite{Str}, and related enquiries from Brian Conrad, for
being reminded of this unneeded gap in the literature and of the
desirability of detailing a straightforward argument.  A much
clumsier version of the story told here is given in
\cite{149}, however with additional introductory details that may be
helpful to the reader.

\begin{theorem}\label{th:2} Set $\mathcal C:
y^2=D(x)$, with
$D$ monic, squarefree, and of degree $2g+2$. Suppose the
divisor at infinity on the Jacobian of the curve
$\mathcal C$ is torsion of order $m\gt1$, equivalently
the domain $\Q[x,y]$ is exceptional in
 containing nontrivial units, and its fundamental unit  $u=a+by$
is of degree
$m$, and say of norm $k$. Denote the
\cfe\ of
$y$ by
$y=\CF{A\\a_1\\a_2\\a_3\\\dots}$. Then, further to
Theorem~\ref{th:1},
\begin{itemize}
\item[(a)] If $\CF{A\\a_1\\a_2\\\dots\\a_{r-1}}=a/b$ with 
$r$ even, then $k=1$.
\item[(b)] If $k=c^2$ is a square, then the \poly\ $b$ factorises
over $\Q$ as say $b=d_+d_-$,  and
$D$ is reducible over $\Q$ because it factorises as the product of
the nontrivial
\poly s
$(a+c)/d^2_+$ and
$(a-c)/d^2_-$. 
\item[(c)] If $k=c^2$ is not a square in $\Q$ then  the \poly\ $b$
factorises over $\Q(c)$ as a product $b=d\conj d$ of \poly s
conjugate over $\Q(c)$, and $D$ factorises over $\Q(c)$ as a product
of the two \poly s $(a+c)/d^2$
and
$(a+\conj c)/\conj d\,^2$. 
\end{itemize} 
\end{theorem}
\noindent For $g=1$ we must have $m=r+1$ by the bounds
\eqref{eq:boundsf}, so the parities of $m$ and $r$ are of
course different; in particular, $m$ odd entails the norm $k=1$. One
readily notices that symmetry implies that always if
$r$ is odd the parities of $m$ and $g$ are different; the
converse is not true if $g\gt1$. For the rest, Theorem~\ref{th:2}
fills in details omitted from Theorem~\ref{th:1}.

An important such `detail', is the observation that if, say, $2\deg
d_+=m$ so $d_+^2=a+c$, then $Dd_-^2=a-c=d_+^2-2c$. So also $d_+
+yd_-$ is a unit of $\Q[x,y]$ plainly contradicting the minimality
of $m$, that is, that $u$ is a fundamental unit.

Furthermore, we see that $D$ has a factor of degree at most $g$ if
the period length $r=2h$ is
even. For then, by conjugation, the line
$$(y+P_h)/Q_h=a_h-(\conj y+P_{h+1})/Q_h$$
is symmetric, that is $P_{h+1}=P_h$, and so $Q_h$ divides $P_h$.
But then $Q_h$ also divides the norm $(y+P_h)(\conj y+P_h)$ and that
entails
$Q_h$ is a factor of $D$.

There are contexts in which one would like to be certain that a
\poly\ $D$ is \emph{not} exceptional. Our results have the
following consequence.
\begin{corollary} If a monic \poly\ $D$ of even
degree is irreducible and with Galois group the
full symmetric group then $D$ is not exceptional;
that is, the \cfe\ of $\sqrt{D}$ is \emph{not} periodic.
\end{corollary}

\end{Remarks*}

\section{Exceptional Polynomials}
\noindent In practice, the start of the \cfe\ of $y=\sqrt D$ quickly
reveals whether or not $D$ is exceptional. For example, it is shown
in
\cite{AR} for $g=1$ that in $y_h=(y+P_h)/Q_h$ the divisor of
$Q_h$ is
$h+1$ times the divisor at infinity. Thus, by well known properties
of Neron--Tate height, the number of decimal digits of the
numerators and denominators of the coefficients of
$Q_h$ (and then also of $P_h$) is $O(h^2)$ unless the divisor at
infinity is torsion. Moreover, in practice that explosion in
complexity of the
$Q_h$ is immediately evident; see \cite{145} for an example.
Moreover, that same explosion in complexity occurs for arbitrary
$g\gt0$ since it follows from addition on the Jacobian of the
curve
$y^2=D(x)$ being given by composition of quadratic forms, that is,
by the
\cfe\ of $y$; 
\cite{Ca} or \cite{La} explain this connection. In any case,
\cite{BCZ}, the matter of explosion of complexity of Pad\'e
approximants of algebraic functions of positive genus is far more
general yet.

In the number field case, the fundamental unit of an order
$\Z[f\omega]$ is some power of the fundamental unit of the
domain of all integers of $\Q(\omega)$. For function fields over a
base field of characteristic zero, however, an order $\Q[x,f(x)y]$
need not possess a unit at all, notwithstanding that $D=y^2$ be
exceptional. In other words, periodicity of $y$ does not at
all guarantee quasi-periodicity of $fy$ for a \poly\ $f$
of positive degree. The requirement in our theorems that $D$  be
squarefree thus really does matter. Specifically, although the
\cfe\ is trivially quasi-periodic for $\deg D=2$, thus $y^2=D$ of
genus
$g=0$, this may not hold for $y^2=f^2D$, even though that curve is
of genus~$0$.  There are interesting papers, see  
\cite{HaHe} and its references, discussing this issue. 

\section{The Quartic Case} \label{g1}
\noindent The case $g=1$ is completely
known over $\Q$, see
\cite{163} and its references, or for example \cite{AZ}. In
particular, one knows by Mazur's Theorem \cite{Ma} that 
the only possibilities for
$m$ are
$m=2$,
$3$ 
$\ldots\,$, $10$, and $12$. From
\cite{158} one learns that in the cases
$m=10$ and $m=12$ it happens that in fact $k=c^2$ never is the
square of a rational. 

For torsion $m\ge4$ one may take $D_m$ as
$(X^2+v-w^2)^2+4v(X+w)$ without loss of generality;
$D_3(X)=(X^2-w^2)^2+4v(X+w)$, while
$D_2(X)=(X^2+u)^2+4w$.

\begin{theorem} \label{th:3} 
Set $\mathcal C_m:
y^2=D_m(x)$, with
$D_m$ monic, squarefree, and of degree $4$. Suppose the
divisor at infinity on the Jacobian of the curve
$\mathcal C_m$ is torsion of exact order $m\gt3$. Then $D_m(x;t)$ is
reducible  over $\Q$ if $m$ is odd or in the cases listed in Table~II. Otherwise, its
Galois group is the dihedral  group~$\mathcal D_4$, other than for the exceptions
listed in Table~I.
\end{theorem}
\begin{proof} We know from above that $D_m(x,t)$ is reducible if $m$ is
odd or if the norm $k_m(t)$ of the fundamental unit happens anyhow to be
a square. Specifically, \cite{158} reports that
$k_8(t)=4(t-1)(2t-1)^2/t^3$, $k_6(t)=4t$, and $k_4(t)=4t$, explaining several of
the entries in Table~II.  Thus we may suppose that $k=c^2$
with
$c$ quadratic irrational over $\Q$.

The Galois group $G_D$ of $D=D_{m}$ is the dihedral group $\mathcal
D_4$ exactly when the zeros of $D$ are $\alpha_1$,
$\alpha_3$,
$\alpha_2=\conj\alpha_3$, and $\alpha_4=\conj\alpha_1$, where
$\conj{\phantom{\alpha}}=(14)(23)$ is conjugation over $\Q(c)$. Then $G_D$ is
generated by that conjuagation  and $\sigma=(1234)$.

Conversely, given that $D$ factorises over $\Q(c)$, the cubic resolvent 
$C_D$ of $D$ must have a rational zero $\alpha_1\alpha_3+\alpha_2\alpha_4$.
The other two zeros
$\alpha_1\alpha_2+\alpha_3\alpha_4$ and $\alpha_1\alpha_4+\alpha_2\alpha_3$
are invariant under the conjugation but are transposed by $\sigma$ and, for
that matter, also by the $4$-cycle $\tau=(1243)$.

If these other zeros of $C_D$ are rational then both $\sigma$ and $\tau$
must be involutions commuting with the conjugation.
Then, recalling that $D$ is irreducible over $\Q$, its Galois group $G_D$
is the Viergruppe $\mathcal V$.
If the pair of zeros is irrational but $D$ factorises over the splitting
field of $C_D$ then $\tau$ generates $G_D$ and the Galois group of $D$  is
the cyclic group
$\mathcal C_4$. Incidentally, we use the helpful remarks 
\cite[Algorithm 4.2 at page~10]{Healy}, explicitly to distinguish the case
$\mathcal C_4$ from $\mathcal D_4$.

\subsection*{Even calculations} We investigate each case $m=12$, $10$, $8$, $6$,
$4$ in detail using the data listed in \cite{158}.
For example, the cases $m=12$ and $m=10$ are
given by 
\begin{multline}\label{eq:fam12}
v_{12}(t)=(t- 1)(2t-1)(3t^2 - 3t + 1)(2t^2 - 2t + 1)/t^4\,;
\\\quad w_{12}(t)=-(6t^4 - 16t^3 + 14t^2 - 6t + 1)/ 2t^3\,;
\end{multline}
\begin{equation}\label{eq:fam10}
v_{10}(t)=\frac{t^3(2t-1)(t-1)}{(t^2-3t+1)^2}\,;\quad
w_{10}(t)=\frac{2t^3-2t^2-2t+1}{2(t^2-3t+1)}\,.\end{equation}
Here the parameter $t$ runs through all `regular'
elements of $\Q$; in both cases the irregular rational values
are $t=1$, $t=1/2$, and $t=0$.

By Theorem~\ref{th:1}(c) we know that $D_m(x;t)$ factorises over
$\Q\(c(t)\)$. If it also factorises over $\Q$ it must do so as a
product $(x^2-px+q)(x^2+px+q')$. One solves (rather, \textsf{Maple}
\cite{MV} solves) this condition for
$p=p(t)$, in each case obtaining two polynomial equations in $p$
and $t$, with one an elliptic curve and the other a quadratic in
an auxiliary variable. The condition that its discriminant be a
square also is an elliptic curve. 

In the case $m=12$, both of these equations ultimately
transform birationally (here \textsf{PARI-GP} \cite{pari} lends a
hand) to the minimal model $y^2=x^3-x^2+x$. This is is \textsf{24A4} in
John Cremona's tables \cite{Cre}; thus with conductor $24$. It has
rank
$0$ and cyclic torsion of order $4$; the torsion points are
$(0,0)$, $(1,1)$, $(1,-1)$, and $\infty$ and correspond to 
irregular values of $t$. So $D_{12}(x;t)$ is irreducible over
$\Q$ for all regular $t\in\Q$.

When, instead, we check the cubic resolvent, for example when
$m=10$, we find that its rational zero is 
$$ {\left (2\,{t}^{3}-4\,{t}^{2}+4\,t-1\right
)\left (2\,{t}^
{3}-4\,{t}^{2}+1\right )}\bigm/{2\left ({t}^{2}-3\,t+1\right )^{2}}
$$
and if the discriminant of the remaining quadratic factor of $C_D$
is a square then the elliptic curve $s^2=(4t^2-2t-1)(2t-1)$ must
have admissible rational points. However, its minimal model
$y^2=x^3+x^2-x$ is \textsf{20A2} in Cremona's tables and it has rank $0$
and cyclic torsion of order $6$. The torsion points are
$(0,0)$, $(\pm1,\pm1)$, and $\infty$
and correspond to irregular values of $t$.

Following \emph{both} the alternative approaches for each of
$m=12$ and
$m=10$ verifies a result we have used above, to wit Tran's result
\cite[p.\,400\emph{ff}]{158} that neither
$\kappa_{12}(t)$ nor $\kappa_{10}(t)$ --- see \S~\ref{ss:quasi}
at page~\pageref{ss:quasi} above --- can be the square of a
rational for regular $t\in\Q$. 

For these and the remaining even cases $m=8$, $m=6$, and $m=4$, where we
know that $k=\kappa_m(t)$ may be a square for some regular $t$, we
followed both approaches and found that when $D_m(x;t)$ is irreducible its
Galois group
$G_D$ is the dihedral group $\mathcal D_4$ except in the cases encapsulated in
the following table.

{\Small\noindent
\begin{center}\begin{tabular}{|c|c|c|c|}
  \hline
   $m$ & $(v,w)$ & $G_D=\mathcal C_2\times \mathcal C_2$ & $G_D=\mathcal
   C_4$\\
   \hline
    & & & \\
  $4$  & $\left(t,\frac12\right)$
       & $t=\frac1{16}({s^2-1})$ & $t=-\frac1{16}/(s^2+1)$  \\
 & & & \\
  $6$  & $\left(t(t-1),1-t/2\right)$
       & $t=8/({9-s^2})$ & $-$ \\
 & & & \\
  $8$  & $\begin{array}{c}\left((t-1)(2t-1),\right.\phantom{abc} 
\\  \phantom{abc}\left.-({2t^2-4t+1})/{2t}\right) \end{array}$ & $-$ & $-$ \\
 & & & \\
  $10$ &
$\begin{array}{c}\left({t^3(2t-1)(t-1)}/{(t^2-3t+1)^2},\right.\phantom{ab} 
\\  \phantom{ab}\left. {2t^3-2t^2-2t+1}/{2(t^2-3t+1)}\right) \end{array}$
& $-$ & $*$ \\
 & & & \\
$12$ & $\begin{matrix}\left({{\left (t-1\right )\left (2\,t-1\right )\left (3\,{t}^{2}-3\,t+
1\right )\left (2\,{t}^{2}-2\,t+1\right )}/{{t}^{4}}},\right.\\
 \phantom{abcdefghijklm}-{
{(6{t}^{4}-16{t}^{3}+14{t}^{2}-6t+1)}/{2{t}^{3}}}\,)\end{matrix}$ & $-$ & $-$  \\
&&&\\
   \hline
\end{tabular}
\end{center}\medskip
}
\centerline{Table~I}

\medskip
\noindent Moreover for $m$ even, $D_m(x,t)$ is irreducible except in the
following cases:
\medskip

{\Small\noindent\begin{center}
\begin{tabular}{|c|c|c|c|c|}
  \hline
   $m$ & $(v,w)$ & $D=f_1f_2$ 
  & $D=f_1f_2f_3$ & $D=f_1f_2f_3f_4$ \\
   \hline
       & & & & \\
  $4$  & $\left(t,\frac12\right)$
       & $t=\begin{cases}-s^2,\\ 4s^4-s^2\end{cases}$ 
       & $t=-\left(\frac{s^2-1}4\right )^2$
       & $t= -\left(\frac{s^3-s}{(s^2+1)^2}\right)^2$
 \\
       & & & & \\
  $6$  & $\left(t(t-1),1-t/2\right)$
       & $t=\begin{cases}1-s^2 \\ \frac{(1+s^2)^2}{3s^2+1}\end{cases}$
       & $t=1-\left(\frac{s^2-1}{s^2+3}\right)^2$
       & $-$ \\
       & & & & \\
  $8$  & $\begin{array}{c}\left((t-1)(2t-1),\right.\phantom{abc} 
\\  \phantom{abc}\left.-({2t^2-4t+1})/{2t}\right) \end{array}$
&
$t=1/({s^2+1})$ & $\dag$ & $-$ \\
       & & & & \\
\hline
\end{tabular}\end{center}\medskip
}
\centerline{Table~II}
\medskip

\noindent The notes $*$ and $\dag$ refer to two special cases we resolved
not to attempt to resolve. We found that rational points $(t,u)$ on the
curve
$$u^2=(t-1)(4t^2-2x-1)(2t-1)(t^2-3t+1)t$$
give rise to cases $D_{10}(x;t)$ with Galois group $\mathcal D_4$; and
rational
points on the curve
$$u^2=(t^4-1)(t^2+2t-1)$$
give cases where $D_{8}(x;t)$ splits into three factors over $\Q$. We expect
that neither curve provides regular rational such $t$.

We leave the degenerate case $m=2$, where
$D(x;u,k)=(X^2+u)^2-k$, as an easy exercise.
\end{proof}

\subsection*{Odd remarks} In the odd cases $m=9$, $m=7$, and $m=5$, the 
final remark following Theorem~\ref{th:2} at page~\pageref{th:2},
together with the detailed \cfe s\footnote{As always, such data must be used
modulo typos. Worse, the notation of \cite{158} is slightly
different from that of here and in \cite{163};  its $v$ is
our $4v$.} in
\cite{158}, shows that 
\begin{gather*}
\(x-\tfrac12({t^3-3t^2+4t-1})\)\quad\text{divides $D_9(x;t)$}\,,
\\
\(x+\tfrac12({t^2-3t+1})\)\quad\text{divides $D_7(x;t)$}\,,\\
\(x-\tfrac12(t+1)\)\quad\text{divides $D_5(x;t)$}\,;
\end{gather*}
here \begin{gather*}
v_9(t)=t^2(t-1)(t^2-t+1),\quad w_9(t)=-\tfrac12({t^3-t^2-1}), 
\qquad \text{$t\in\Q\setminus\{0,1\}$}\,,\\
v_7(t)=t^2(t-1),\quad w_7(t)=-\tfrac12(t^2-t-1), 
\qquad \text{$t\in\Q\setminus\{0,1\}$}\,,\\
v_5(t)=t,\quad w_5(t)=-\tfrac12(t-1), \qquad
\text{$t\in\Q\setminus\{0\}$}\,.
\end{gather*}
For completeness we remark that in these cases the residual cubic factor
$G_m(x;t)$ is reducible in the case $m=5$ and $t=s^2(s+1)/(s+1)$ and that then
the surviving quadratic factor is irreducible. With finitely many possible
exceptions,  namely
unlikely rational points on certain curves\footnote{
The respective discriminants $F_m(t)$ of the cubic factors
are  $F_7(t)=t(t-1)(t^3-8t^2+5t+1)$,
$F_9(t)=t(t-1)(t^2-t+1)(t^3-6t^2+3t+1)$, and $F_5(t)=t(t-1)(t^3-8t^2+5t+1)$.
The last case is Cremona's curve \textsf{20A2}, which has rank~$0$ and
torsion~$2$. We saw that $G_7(x;t)$ is irreducible because a putative rational
zero corresponds to a rational point on the curve \textsf{14A4} with rank~$0$
and torsion~$2$. We found a complicated genus~$2$ curve not warranting
report whose rational points might allow $G_9(x;t)$ to factorise.%
} 
of genus more than~$1$, the Galois
groups of the irreducible
$G_m(x;t)$ is always
$\mathcal S_3$.

 The case
$m=3$ is degenerate; however, plainly
$$
D_3(x;v,w)=(x^2-w^2)^2+4v(x+w)=(x+w)(x^3-wx^2-w^2x-4v+w^3)=:(x+w)F\,.
$$ 
If $F$ is irreducible, then its Galois group is $\mathcal A_3$ if
and only if $v=8t^2w^3/(27t^2+1)$.
Further, $F$ has a zero $r$ when
$v=(w+r)(w-r)^2/4$; specifically
$$F= (x-r)(x^2-(w-r)x-w^2-rw+r^2).$$
$F$ splits as the product of three linear factors when
$v= 8w^3(s^2-1)^2/\((s^2+3)^3\)$.
The reader may find it a useful
exercise to extract other details.

\bibliographystyle{amsalpha}


\label{page:lastpage}
\end{document}